\documentclass[11pt]{article}

\usepackage{graphicx}
\usepackage{graphicx,psfrag}
\usepackage{amsmath,amssymb,amsthm}
 
\def \G {\Gamma}

\def \n {\noindent}

\def \E {{\bf  \, E  }}

\def \Tr {{\mathrm{\, Tr}}}

\def \BE {{\bf E}}
\def \bC {{\mathbb C}}

\def \bE {{\mathbb E}}

\def \bN {{\mathbb N}}

\def \bR {{\mathbb R}}

\def \CA {{\cal A}}
\def \CB{{\cal B}}

\def \CC {{\cal C}}

\def \CH {{\cal H}}
\def \CU {{\cal U}}

\def \CM{{\cal M}}
\def \CW {{\cal W}}

\def \CP {{\cal P}}

\def \CV {{\cal V}}

\def \D {\Delta}

\def \G{{\Gamma}}

\def \a {\alpha}

\def \b {\beta}

\def \s {\sigma}

\def \vep {\varepsilon}

\def \ra {\rangle}

\def \ti {\tilde}

\def \RH {{\mathrm H}}

\def \RI {{\mathrm I}}

\def \RA {{\mathrm A}}

\def \RB {{\mathrm B}}

\def \rb {{\mathrm b}}

\def \ra {{\mathrm a}}

\def \rb {{\mathrm b}}

\def \det {\hbox{det}}

\begin{document}
\title{On asymptotic properties of Bell polynomials and 
concentration of vertex degree of large random graphs
}

\author{O. Khorunzhiy\\ Universit\'e de Versailles - Saint-Quentin \\45, Avenue des Etats-Unis, 78035 Versailles, FRANCE\\
{\it e-mail:} oleksiy.khorunzhiy@uvsq.fr}
%\date{}
\maketitle
\begin{abstract}
We study  concentration properties of vertex degrees of  $n$-dimensional 
Erd\H os-R\'enyi random graphs with the edge probability $\rho/n$
%in the limit  when $n$ and $\rho$ tend to infinity  
by means of high moments 
of these random variables in the limit when $n$ and $\rho$ tend to infinity. These moments
 are asymptotically close to  one-variable Bell polynomials 
$\CB_k(\rho), k\in \bN$ that represent  moments of the Poisson probability distribution 
$\CP(\rho)$. 

We study  asymptotic  behavior of the Bell polynomials 
and modified Bell polynomials
for large  values of $k$ and $\rho$
with the help of the local limit theorem for  auxiliary
random variables.
 Using the results 
obtained, 
  we get   the upper bounds for  the deviation probabilities 
of the normalized maximal vertex degree of the  Erd\H os-R\'enyi random graphs 
 in the limit $n,\rho\to\infty$  such that the ratio $\rho/\log  n $ remains finite or infinitely increases. 
  
\end{abstract}

%{\it Running title:}  On Ihara zeta function of random graphs

%%%%%%%%%%%%%%%%%%%%%%%%%%%%%%%%

\section{Introduction}

This paper is motivated by  studies  of the spectral properties of random matrices 
associated with random graphs of the Erd\H os-R\'enyi type \cite{B}. 
In these graphs the edges can be 
 represented by a family of  independent
 Bernoulli random variables and 
we consider the case when  the average value of these variables
is given by $\rho/n$, where $n$ is the number of vertices in the graph 
$\G_n^{(\rho)}$.

More precisely, we assume that the edges of $\G_n^{(\rho)}$ are non-oriented and there are no loops in $\G_n^{(\rho)}$. 
In this case the random graph $\G_n^{(\rho)}$ is associated with 
its adjacency matrix that is  an $n$-dimensional real symmetric matrix $\RA_n^{(\rho)}$, whose elements 
 above the 
diagonal are given by an ensemble of jointly independent 
Bernoulli random variables $\{\ra_{ij}^{(n,\rho)}\}_{1\le i< j\le n}$
$$
\left( \RA_n^{(\rho)}\right)_{ij} = \ra_{ij}^{(n,\rho)}  =
 \begin{cases}
  1 , & \text{with probability $\rho/n$} , \\
0, & \text {with probability $1-{\rho/n}\, $} 
\end{cases}, \   1\le i<j\le n, \  0<\rho<n
\eqno (1.1)
$$ 
and $\ra_{ii}^{(n,\rho)}=0, 1\le i\le n$.   
In the majority of aspects, the random graphs ensemble   $\{\Gamma_n^{(\rho)}\}$ is very similar  to the  random graphs
introduced and studied by P. Erd\H os and A. R\'enyi \cite{ERRG} and  we refer to $\{\Gamma_n^{(\rho)}\}$  as to the Erd\H os-R\'enyi
  ensemble of random graphs   \cite{B,JLR}.

Asymptotic properties of random graphs  $\G_n^{(\rho)}$ in the limit of infinite $n$ are extensively studied (see monographs \cite{Durr,JLR,Pal}), as well as the spectral properties of  their adjacency matrices $\RA_n^{(\rho)}$ \cite{FK,J,KS,McK}. 
The spectral properties of the second differential form on the vertices of $\G_n^{(\rho)}$ given by 
$$
\Delta_{\G_n^{(\rho)}} = \RB_n^{(\rho)} - \RA_n^{(\rho)},
$$
where $\left(\RB_n^{(\rho)}\right)_{ij}  = \delta_{ij} \rb_i ^{(n,\rho)}$, 
$1\le i\le j\le n$,
$$
\rb_i^{(n,\rho)} = \sum_{l=1}^n \ra_{il}^{(n,\rho)}
\eqno (1.2)
$$
and $\delta_{ij}$ is the Kronecker delta-symbol
 equal to one if $i=j$ and to zero otherwise have been also extensively studied  
 \cite{DJ,KSV}. 
 Usually  $\Delta_{\G_n^{(\rho)}}$ 
 is referred to as a discrete version of the Laplace operator determined 
 on  graph $\G_n^{(\rho)}$. Let us note that random variables 
 $\{\rb_i^{(n\rho)}\}_{i=1,\dots, n}$ (1.2) represent  degrees of vertices of $\G_n^{(\rho)}$ and 
 obviously follow the  binomial distribution ${\cal B}(n-1,\rho/n)$.

Normalized versions of (1.2) are  important, in particular,   
because of their role in the Ihara zeta function of $\G_n^{(\rho)}$. This function  
 can be defined, in its determinant form,
by  the following relation \cite{Ba,I,ST}      
$$
Z_{\Gamma_n^{(\rho)}}(u)
= \left( (1-u^2)^{r-1} \ \det\,  \RH_n^{(\rho)}(u)\right)^{-1},
\quad u\in \bC, 
\eqno (1.3)
$$
where 
$r-1= \Tr (\RB_n^{(\rho)}-2{\hbox{I}})/2, \ (\RI)_{ij} = \delta_{ij}$,
and 
$$
\RH_n^{(\rho)}(u)=  u^2\RB_n^{(\rho)}- u\RA_n^{(\rho)} + (1-u^2){\RI}.
\eqno (1.4)
$$ 
The Ihara zeta function  is 
a kind of exponential generating function of  walks over the graph 
$\G_n^{(\rho)}$ \cite{HST} and  in general, is rather 
difficult to be computed. 
The random matrix approach can be useful
in  the studies of the Ihara zeta functions and related properties of graphs
either in average or with probability one \cite{K-17}.

Logarithm of  determinant 
of (1.3)  
can be expressed in terms  of  the normalized eigenvalue counting function 
of  matrices $ \RH_n^{(\rho)}(u)$ (1.4).
In paper \cite{K-17}, it is shown 
that this   eigenvalue distribution  function of matrices 
$\RH_n^{(\rho)}(u)$ of the Erd\H os-R\'enyi random graphs 
$\{ \Gamma_n^{(\rho)}\}$, regarded with  the normalization parameter $u= v/\sqrt \rho$,
converges as  
  $n, \rho\to\infty $ to a measure 
  given by a shift of the Wigner semi-circle law \cite{W}. 
  This limiting measure determines the limiting expression for 
  the logarithm of the Ihara zeta function (1.3), if it exists.
  To study  convergence of 
 the logarithm of Ihara zeta function 
 $\log Z_{\Gamma_n^{(\rho)}}(v/\sqrt \rho)$ 
  by itself, one should 
   specify the convergence of  the eigenvalue distribution function of (1.4) 
   to the limiting measure and, in particular, to  study how rapidly 
  the 
   diagonal elements of $\RH^{(\rho)}_n(v/\sqrt \rho)$  given by 
  $
   {v^2  \rb_i^{(n,\rho)}/ \rho}
  $ 
  \mbox{converge to $1$}
   as $n,\rho\to\infty$. 
  This question  is important also in the studies of the spectral properties of random adjacency matrices 
  $\RA_n^{(\rho)}$  \cite{KS}.

   To obtain upper bounds for the deviation probabilities of centered random variables 
   $\tilde  \rb_i^{(n,\rho)} =  \rb_i^{(n,\rho)} - (n-1)\rho/n$, we consider
    their moments that are  close to the  modified Bell polynomials $\tilde \CB_k(\rho)$ in the limit $n\to\infty$.
  The main technical results of the present article
  are related with  
 asymptotic properties of 
  polynomials $\tilde \CB_k(\rho)$ in the limit when $k$ and $\rho$ both tend to infinity.
  These results are new and can  have their proper interest because
  the family $\{\tilde \CB_k(\rho)\}_{k\in \bN} $ represents the moments of the centered Poisson distribution $\tilde \CP(\rho)$ \cite{KPST}
  whose asymptotic behavior, up to our knowledge, has not yet been  studied. 
  
  Using the results on asymptotic behavior of centered Poisson distribution, 
  we obtain, with the help of the Markov-type inequalities, upper bounds 
 on the deviation probabilities of the maximal vertex  degree of large Erd\H os-R\'enyi random graphs. Corresponding concentration results
 represent the main subject of the present article.

The paper is organized as follows. In Section 2, we give rigorous definitions 
of  binomial and Poisson random variables and the Bell polynomials
that we consider. Then we formulate our statements concerning    asymptotic
properties of the Bell polynomials $\CB_k(x)$  and modified Bell polynomials
$\tilde \CB_k(x)$ in the limit $k,x\to\infty$. These statements represent our 
main technical results. In Section 3,
we use them to study  concentration properties of vertex degrees
of large random graphs. We compare the  results obtained with known facts, in particular with the Erd\H os-R\'enyi law of large numbers 
 \cite{ER}.

In \mbox{Section 4,} we prove theorems
formulated in Section 2. 
We introduce auxiliary random variables $Z$ and $\tilde Z$ related with 
$\CB_k(x)$  and 
$\tilde \CB_k(x)$, respectively, and prove the local limit theorems for these random variables that give the limiting expressions  
for $\CB_k(x)$  and 
$\tilde \CB_k(x)$. 
Section 5 contains the proofs of auxiliary statements.

%%%%%%%%%%%%%%%%%%%%%%%%%%%%%%%%%%%%%%%%%%%%%%%%%%%%%%%%%%%%%%%%%%%%%%

\section{Bell and restricted Bell polynomials }

The material of this section can be of
a general interest, so we replace the vertex degrees (1.2) by a 
 sum of $n$  jointly independent Bernoulli random variables 
$$
a^{(n,\rho)}_j =
 \begin{cases}
  1 , & \text{with probability $\rho/n$} , \\
0, & \text {with probability $1-{\rho/n}\, $} 
\end{cases}, \quad 1\le j\le n, \  \rho>0
\eqno (2.1)
$$
given by 
$$
X_n^{(\rho)} = \sum_{j=1}^n a^{(n,\rho)}_j.
$$
We also consider the centered random variables 
$$
\tilde X_n^{(\rho)}= X_n^{(\rho)} - \bE X_n^{(\rho)} = 
\sum_{j=1}^n \tilde a^{(n,\rho)}_j= \sum_{j=1}^n \left(a^{(n,\rho)}_j - \rho/n\right).
$$
Here and below we denote by $\bE $ the mathematical expectation with respect to the measure generated
by the family $\CA^{(n,\rho)}= \left\{ \{a_j^{(n,\rho)}\}_{1\le j\le n}\right\} $.

 Probability distributions of  $X_n^{(\rho)}$ and $\tilde X_n^{(\rho)}$ 
 converge as $n\to\infty$ to the Poisson  distribution  $\CP(\rho) $ 
 and to the centered Poisson distribution    $\tilde \CP(\rho)$, respectively
 (see Section 5 for details).
 It is known that the moments of $\CP(\rho)$ and $\tilde \CP(\rho)$ of order $k$
are given by the one-variable Bell polynomials $\CB_k(\rho)$ \cite{Car}
and their modified versions $\tilde \CB_k(\rho)$ \cite{KPST,P}. 
Therefore the $k$-th moments of $X_n^{(\rho)}$ and $\tilde X_n^{(\rho)}$ 
converge to $\CB_k(\rho)$ 
and  $\tilde \CB_k(\rho)$ as $n\to\infty$. 
We will need the following generalization of this observation.

\vskip 0.2cm 
\n  {\bf Lemma 2.1.} 
 {\it
 Let 
 $$
 \CM_k^{(n,\rho)}=  \bE \left( X_n^{(\rho)}\right)^k \quad {\hbox{and}} 
 \quad 
\tilde \CM_k^{(n,\rho)}=\bE \left(\tilde  X_n^{(\rho)}\right)^k.
\eqno (2.2)
$$
If $n$ infinitely increases, then 
$$
\CM_k^{(n,\rho)} =  {\cal B}_k(\rho) (1+o(1)), \quad n\to\infty
\eqno (2.3)
$$
for any finite or infinite $\rho$ and $k= o(\sqrt n)$,  
and 
$$
\tilde \CM_k^{(n,\rho)} =  \tilde {\cal B}_k(\rho) (1+o(1)), \quad 
n\to\infty
\eqno (2.4)
$$
for any sequence $k, \rho$ such that $k= o(\sqrt n)$ and $k\rho = o(n)$, 
where 
$$
{\cal B}_k(x) = \sum_{(l_1,l_2,\dots,l_k)}^k B^{(k)}_{l_1, l_2 ,\dots, l_k}\, x^{l_1+ l_2+\dots +l_k},  \quad k\ge 1, \ x\in \bR,
\eqno (2.5)
$$
with
$$
B^{(k)}_{l_1,l_2,\dots, l_k} = {k!\over (1!)^{l_1} l_1! \, (2!)^{l_2} l_2! \cdots (k!)^{l_k}l_k!}
\eqno (2.6)
$$
and 
$$
\tilde {\cal B}_k(x) = \sum_{(l_2,l_3, \dots,l_k)'} \tilde  
B^{(k)}_{ l_2 ,\dots, l_k}\,  x^{l_2+l_3+\dots +l_k},
\quad k\ge 1, \ x\in \bR,
\eqno (2.7)
$$ 
with 
$$
\tilde  B^{(k)}_{l_2,l_3,\dots, l_k} = 
{k!\over (2!)^{l_2} l_2! (3!)^{l_3} l_3! \cdots (k!)^{l_k}l_k!}.
\eqno (2.8)
$$
The sum in (2.5) 
runs over such integers $l_i\ge 0$
that $l_1+ 2l_2 + \dots +kl_k = k$ and 
the sum in (2.7) runs over such $l_i\ge 0$  
that $2l_2+3l_3+\dots + kl_k=k$.
}

%We refer to 
%$\tilde \CB_k(\rho)$
%as to the Bell polynomials of restricted Bell numbers, 
%sor simply as to restricted Bell polynomials.

\vskip 0.2cm 

We prove Lemma 2.1 in Section 5 below. 
The value $\CB_k(1)= B_k$ represents the number of all possible partitions
of the set of $k$ elements into non-empty subsets (blocks) and 
the numbers $B_k,{k\ge 0}$
are known as  the  Bell  numbers
  \cite{B1,B2} 
 (see also  \cite{Berndt,Touch}).
The value $\tilde B_k= \tilde \CB_k(1)$ represents the number of all possible partitions
of the set of $k$ elements into blocks of cardinality strictly greater than one. 
On can say  that  $\tilde  B_k$ are 
 restricted versions of the Bell numbers.
We  refer to   $\tilde {\cal B}_k(x)$ as to  the Bell polynomials of restricted  partitions,  or 
 simply as to restricted Bell polynomials.

 \vskip 0.2cm 
As we have said above, we are interested in the asymptotic behavior of 
polynomials
 $\tilde {\cal B}_k(x), x>0$
 in the limit 
when $k$  and $x$ both tend to infinity. We use a method 
that can be applied to the Bell polynomials and  to their modified versions as well.
It would be interesting to   compare
results obtained in these two cases, so we formulate them in parallel.
Let us start with  the classical one-variable Bell polynomials $\CB_k(x)$.

%%%%%%%%%%%%%%%%%%%%%%%%%%%%%%%%%%%%%%%%%

 \vskip 0.2cm 
\noindent {\bf  Theorem 2.1.} {\it  The Bell polynomials $\CB_k(x)$, $x>0$ have the following properties:

\vskip 0.2cm 
\noindent a) if $x$ is non-vanishing as $k\to\infty$ and $x= o( k)$,
 then 
$$
\CB_k(x) = \left( {k\over e (\ln k - \ln x)} (1+o(1))\right)^k, \quad k\to\infty;
\eqno (2.9)
$$
\vskip 0.2cm 
\noindent  b) if sequence $(x_k)_{k\in \bN}$ is such that $x_k/k \to\chi>0$  as $k\to\infty$, then 
$$
{\cal B}_k(x_k) = \big(k\chi e^v \left(1+o(1)\right)\big)^k , \quad k\to\infty,
\eqno (2.10)
$$
where 
$$
v=h(u)=  u - 1 + {1\over u} - {1\over u e^u}
\eqno (2.11)
$$ 
and $u= u(\chi) $ verifies equality   
$$u e^u={1\over  \chi}
\eqno (2.12) 
$$ known as the Lambert equation \cite{DB};
moreover, if $x_k = \chi k$, $k\in \bN$, then 
$$
{\cal B}_k(x_k) = \big(k\chi e^v \big)^k (1+o(1)), \quad k\to\infty;
\eqno (2.13)
$$

\vskip 0.2cm 
\noindent  c) if sequence $(x_k)_{k\in \bN}$ is such that  and $ x_k/k = \chi_k \to\infty$,  then 
$$
{\cal B}_k(x_k) = 
 \left({k \chi_k}  \big( 1+o(1)\big) \right)^k, \quad  \ k\to\infty.
 \eqno (2.14) 
$$}

{\it Remark.} Relation (2.9) remains  true also in the asymptotic regime when $x\to 0$ at the same time as 
$k\to\infty$ under certain additional condition; for example,
we can show this  in the case when $x\gg k e^{-k^{1/7}}$ as $k\to\infty$,
 see subsection 4.1 below.

\vskip 0.2cm

Asymptotic properties of 
Bell polynomials  $\CB_k(x)$, $k\to\infty$ with negative $x<0$    
    have been studied earlier
\cite{D,E,Zh} with the help of the  differential equations techniques.
Here we develop a method of \cite{TE} based on the local limit theorems for an auxiliary random variable.
It should be mentioned that paper \cite{E} reports  also results on  the asymptotics of 
$\CB_k(kz)$
as $k\to\infty$ for $z$ belonging to a compact subset of $\bC\setminus [-e,0]$
thus matching  
 our expressions (2.9), (2.10) and (2.13).

Regarding
 restricted Bell polynomials $\tilde \CB_k(x)$, we observe that 
 asymptotic properties  of $\tilde \CB_k(x)$ 
 are similar to those of $\CB_k(x)$  in 
the limiting transitions  (a) and (b) of Theorem 2.1. In contrast, 
 in the third and the most important for us limiting transition when 
$x_k/k\to\infty$ as $k\to\infty$, the asymptotic behavior of $\tilde \CB_k(x_k)$ essentially differs
from that of $\CB_k(x_k)$ (2.12). Let us formulate our main technical results.

\hskip 0.2cm

\noindent {\bf Theorem 2.2.}  
{\it Restricted Bell polynomials  $\tilde \CB_k(x)$, $x>0$  have the following properties:

\vskip 0.2cm 
\noindent a) if $x$ is non-vanishing as $k\to\infty$ and  $x = o( k)$,
 then 
$$
\tilde \CB_k(x) = \left( {k\over e (\ln k - \ln x)} (1+o(1))\right)^k;
\eqno  (2.15)
$$

\vskip 0.2cm

\noindent  b) if sequence $(x_k)_{k\in \bN}$ is such that  $x_k/k\to \chi>0$ when $k$ tends to infinity, then 
$$
\tilde {\cal B}_k(x) = \left({k\chi e^{\tilde v} (1+o(1))}\right)^k, 
\quad k\to\infty,
\eqno (2.16)
$$
where 
$$
\tilde v=  \tilde h (u) = u - 1 + {1\over u} + \ln {e^u-1\over  u}
 \eqno (2.17) 
  $$ 
and $u=u(\chi)$ is determined by the following Lambert-type equation,
$$
u(e^u-1) = {1\over  \chi};
\eqno (2.18)
$$
moreover, if $x_k = \chi k$, $k\in \bN$, then 
$$
\tilde {\cal B}_k(x) = \big({k\chi e^{\tilde v} }\big)^k (1+o(1)), 
\quad k\to\infty,
\eqno (2.19)
$$

\vskip 0.2cm 
\noindent  c) if sequence $(x_k)_{k\in \bN}$ is such that 
 $\chi_k = x_k/k\to\infty$ and $x_k /k^{3/2}\to 0$ as $k\to\infty$,  then 
$$
\tilde {\cal B}_k(x) = 
 \left(k \sqrt{ {\chi_k\over  e}}  \big( 1+o(1)\big) \right)^k
, 
 \quad k\to\infty.
 \eqno (2.20)
$$}

{\it Remark.} Relation (2.15) remains true in the asymptotic regime 
when $x$ vanishes at the same time as $x\to\infty$ under certain additional condition; in particular, we can show this in the case when 
$x\gg k e^{-k^{1/15}}$ as $k\to\infty$, see subsection 4.3 below. 

\vskip 0.2cm 
The difference between the right-hand sides of (2.14) and (2.20) can be easily explained by the fact that 
in the limit when $x/k\to\infty, k\to\infty$,
the leading contributions to $\CB_k(x)$ and $\tilde \CB_k(x)$
are given by  terms with the maximal degree of $x$. In the  case
of Bell polynomials this is simply $x^k = (k \chi_k)^k$ while in the  case
of restricted Bell polynomials, this term is determined by, roughly speaking, 
$$
x^{k/2} {k!\over 2^{k/2} (k/2)!} = \left({xk\over e} 
\big(1+o(1)\big)\right)^{k/2}.
\eqno (2.21)
$$
When deriving (2.21), we have taken into account that the number of 
partitions of $2m$ elements into $m$ blocks of cardinality 2 is given by 
$(2m)!/(2^m m!)$ and  then  used the Stirling formula
$$
k! = \sqrt{2\pi k}\left( {k\over e}\right)^k \big(1+o(1)\big).
\eqno (2.22)
$$
The right-hand side of (2.21) is obviously in accordance with that of  (2.20).
We prove Theorems 2.1 and 2.2  in Section 4 below. 

Let us note that relations (2.10) and (2.16) can be rewritten in the following forms,
$$
\lim_{k\to\infty,\  x/k\to\chi} \ {1\over k} \ln 
\left( {\CB_k(x)\over x^k}\right) = \psi(\chi)
\eqno (2.23)
$$
and 
$$
\lim_{k\to\infty,\  x/k\to\chi} \ {1\over k} \ln 
\left( {\tilde \CB_k(x)\over x^k}\right) = \tilde \psi(\chi),
\eqno (2.24)
$$
where $\psi(\chi) = {h(u)}$ with $h(u)$ determined by (2.11) and (2.12)
and $\tilde \psi(\chi) = {\tilde h(u)}$ with 
$\tilde h(u)$ determined by (2.17) and (2.18).

%%%%%%%%%%%%%%%%%%%%%%%%%%%%%%%%%%%%%%%%%%%%%%%%%%%%%%%%%%%%%%%
%%%%%%%%%%%%%%%%%%%%%%%%%%%%%%%%%%%%%%%%%%%%%%%%%%%%%%%%%%%%%%%

\section{Asymptotic properties  of  vertex degree}

Let us return to the 
 random graphs $\Gamma_n^{(\rho)}$ whose adjacency matrix 
$
\RA_n^{(\rho)}
$
is given by (1.1)
and the elements $\left(\RA_n^{(\rho)}\right)_{ij}$ with $1\le j<i\le n$ are determined 
 by the symmetry condition. 
 Random variables $\{\ra_{ij}^{(n,\rho)}\}_{ 1\le i< j\le n}$ are jointly independent and we consider an infinite triangle array 
 $\CA'=\left\{  \{ \ra_{ij}^{(n,\rho)}\}_{ 1\le i\le j\le n}, n\in \bN
 \right\}$
assuming that 
 the values of $\rho= \rho_n$ are determined for any $n\in \bN$.  
 We denote the mathematical expectation
 with respect to the measure generated by $\CA'$ by $\bE'$. 

Diagonal elements $\rb_i^{(n,\rho)}$ of  matrix $\RB_n^{(\rho)}$ (1.3)  represent degrees  
 of vertices $v_i$ of random graph $\Gamma_n^{(\rho)}$.
 The aim of this section is to 
study the convergence of    the maximal vertex degree
$$
  d^{(n,\rho)}_{\max} = \max_{1\le i\le n} \{ \rb_i^{(n,\rho)}\}
  \eqno (3.1)
$$
to its mean value in the limit of infinite  $n$ and 
$\rho$.
To do this, we consider the centered random variables
$$
\tilde \rb_i^{(n,\rho)} = \rb_i^{(n,\rho)} - \rho(n-1)/n
\qquad {\hbox{and} } \qquad 
\tilde d^{(n,\rho)}_{\max} = d^{(n,\rho)}_{\max}  - \rho(n-1)/n
  \eqno (3.2)
$$
and write the classical 
inequality of the Markov type,
$$
P\left( \vert \tilde \rb_i^{(n,\rho)}\vert > s\right)\le s^{-2k}\  
\bE'\left( \tilde \rb_i^{(n,\rho)}\right)^{2k}. 
\eqno (3.3)
$$ 
It is easy to see that $\tilde \rb_i^{(n,\rho)}$ and 
$\tilde \rb_{i'}^{(n,\rho)}$ are identically distributed and 
that their  moments $\tilde M_k^{(n,\rho)}$ can be written as follows,  
$$
\tilde M_k^{(n,\rho)} 
= \bE' \left(\tilde \rb_1^{(n,\rho)} \right)^k
= \bE \left( \sum_{j=1}^{n-1} \tilde a_j^{(n,\rho)}\right)^k,
 \eqno (3.4)
$$
where $\tilde a_j^{(n,\rho)}= a_j^{(n,\rho)}- \rho/n$ and the random variables $a_j^{(n,\rho)}$ are determined  by  (2.1).
The following statement is analogous to Lemma 2.1. 
 \vskip 0.2cm
{\bf Lemma 3.1.}
{\it If $k= o(\sqrt n)$ and $k\rho = o(n)$ as $n\to\infty$, then 
$$
\tilde M_k^{(n,\rho)} = \tilde \CB_k(\rho)(1+o(1)), \quad n\to\infty.
\eqno (3.5)
$$}
We prove Lemma 3.1 in Section 5 in parallel with the proof of Lemma 2.1.
Relation (3.3) together with  (3.5) 
implies that if $k= o(\sqrt n)$ and $k\rho = o(n)$, then 
$$
P(\vert \tilde \rb^{(n,\rho)}_1\vert > s) \le 
 s^{-2k} \, \tilde \CB_{2k}(\rho) (1+o(1)), \quad n\to\infty.
\eqno (3.6)
$$
Now we can 
use  results of Theorem 2.2 and determine
upper bounds for the deviation probabilities of random variables (3.2).
%$\tilde \rb_1^{(n,\rho)} $.
While our main interest is related with the maximal 
vertex degree $d_{\max}^{(n,\rho)}$ (3.1), we start with a detailed  account for the
random variable $\rb_1^{(n,\rho)}$ whose probability distribution
is  the  binomial one ${\cal B}(n-1,\rho/n)$.

\vskip 0.2cm 

{\bf Theorem 3.1.}
{\it If
and $\rho_n = \tilde \chi \ln n$ as $n\to\infty$, then 
$$
\lim_{n\to\infty} P\left(  
\vert {\rb_1^{(n,\rho)}/ \rho_n} - 1\vert \ge   \ti s \right) 
= 0
\eqno (3.7)
$$ 
for any $ \ti s> e^{\ti v}$, where $\ti v = \ti h (u)$ is given by (2.17)
and $u$ is determined by equation (2.18),
$$
u(e^u - 1) = {1\over \tilde \chi}.
$$
For any $ \ti t> e^{\ti v +1}$, the following relation is true, 
$$
P\left(\limsup_{n\to\infty} \left\{\omega:\  
\vert {\rb_1^{(n,\rho)}/ \rho_n} - 1\vert \ge   \ti t \right\}\right) 
= 0.
\eqno (3.8)
$$

\vskip 0.2cm 
Proof.} Let us consider   (3.6) 
with $s= s' \rho_n$, $\rho= \rho_n$ and $k= \lfloor \ln n\rfloor$. Then obviously 
$\rho_n/k\to\tilde \chi$, relation (2.16) takes the form
$$
\tilde \CB_{2k} (\rho_n) = \left( \rho_n e^{\ti v}\big(1+o(1)\big)\right)^{2k}
$$
and we deduce from (3.6) that 
$$
P(\vert \tilde \rb^{(n,\rho)}_1\vert > \tilde s' \rho_n) \le \left( 
{e^{\ti v}\over \ti s'}\big(1+o(1)\big)\right)^{2\lfloor \ln n\rfloor }= \exp\left\{ - 2\lfloor \ln n\rfloor 
\big(\ln \ti s' - \ti v\big)\right\}.
$$
Let $\delta $ be such that $\ti s = e^{\ti v} (1+2\delta)$. Taking 
 $\ti s' = e^{\ti v}(1+\delta)$, we get  inequality
 $$
 P(\vert \tilde \rb^{(n,\rho)}_1\vert > \tilde s' \rho_n)
\le \exp\left\{- 2 \lfloor \ln n\rfloor \ln (1+\delta) \right\} = 
n^{-2 (\lfloor \ln n\rfloor / \ln n)\ln (1+\delta)}.
\eqno (3.9)
$$
It is clear that 
$$
P\left(  
\vert {\rb_1^{(n,\rho)}/ \rho_n} - 1\vert \ge   \ti s \right) 
\le P\left(  
\vert {\rb_1^{(n,\rho)}/ \rho_n} - (n-1)/n\vert \ge   \ti s' \right) 
= P\left(\vert \ti \rb_1^{(n,\rho)}\vert \ge \ti s'\rho\right)
$$
for all $n\ge n_0$, where $n_0$ is such that
$n_0\ge \max\{3, (\delta  e^{\ti v})^{-1}\}$.
This observation together with the upper bound  (3.9) 
implies (3.7).

To prove the second part of Theorem 3.1 given by (3.8), 
we consider $ \delta'$  such that  
$\ti t = e^{\ti v +1}(1+2\delta')$. If $\ti t' = e^{\ti v +1} (1+\delta')$, 
then  
$$
P\left(  
\vert {\rb_1^{(n,\rho)}/ \rho_n} - 1\vert \ge   \ti t \right) 
\le P\left(  
\vert {\rb_1^{(n,\rho)}/ \rho_n} - (n-1)/n\vert \ge   \ti t' \right) 
$$
for all $n\ge n_0$, $n_0 \ge \max\{3, (\delta'  e^{\ti v+1})^{-1}\}$.
Then it follows from (3.9) that
$$
P\left(  
\vert {\rb_1^{(n,\rho)}/ \rho_n} - (n-1)/n\vert \ge   \ti t' \right) 
\le n^{-2 (\lfloor \ln n\rfloor / \ln n)(1+ \ln (1+\delta'))}
$$
and therefore  
$$
\sum_{n= n_0}^\infty  P\left(  
\vert {\rb_1^{(n,\rho)}/ \rho_n} - 1\vert \ge   \ti t \right)  < \infty,
\quad n_0 \ge \max\{3, (\delta'  e^{\ti v+1})^{-1}\}.
$$
The Borel-Cantelli lemma implies (3.8). Theorem 3.1 is proved. $\Box$
\vskip 0.2cm 
Let us note that Theorem 3.1 can be reformulated in terms of 
classical sums of Bernoulli random variables $a_i^{(n,\rho)}$ (2.1) and relation (3.8)
can be rewritten as follows,
$$
P\left(\limsup_{n\to\infty} \left\{\omega:\  
\vert \sum_{i=1}^n a_i^{(n,\rho_n)} - \rho_n \vert \ge   \ti t \rho_n\right\}\right) 
= 0, \quad \rho_n = \ti \chi \ln n,
$$ 
where $\ti t $ is determined by the same condition as in Theorem 3.1.
\vskip 0.2cm

\vskip 0.2cm 
{\bf Theorem 3.2.} {\it If $\rho_n = \chi_n \ln n$ and $\chi_n\to\infty$
as $n\to\infty$, then 
$$
P\left( \lim_{n\to\infty} {\rb_1^{(n,\rho)}\over \rho } = 1\right) = 1.
\eqno (3.10)
$$

Proof.} Since $\rho_n/\ln n = \chi_n \to\infty$, we can use the third result of Theorem 2.2. 
 Relations (3.3) and (3.6) together with (2.20) imply that 
 for any positive $\mu$
$$
P\left(\vert \ti \rb_1^{(n,\rho)}/\rho\vert \ge \mu\right)=
% P\left( \vert \rb_1^{(n,\rho)}/ \rho - (n-1)/n \vert \ge \mu\right)\le 
 \left({1+o(1)\over \mu \sqrt \chi_n}\right)^{2k}
\le {1\over n^{2(\lfloor \ln n\rfloor / \ln n)  \ln(\mu (1+o(1))\sqrt{ \chi_n})}},
\quad n\to\infty,
\eqno (3.11)
$$
where we have omitted the subscript in $\rho_n$ and set  $k= \lfloor \ln n \rfloor$.
Then there exists $n_1$ such that
$$
\sum_{n=n_1}^\infty 
P\left(\vert \ti \rb_1^{(n,\rho)}/\rho\vert \ge \mu\right)< \infty.
$$
The Borel-Cantelli lemma implies equality
$$
P\left( \limsup_{n\to\infty} 
\left\{ \omega: \ 
\vert \ti \rb_1^{(n,\rho)} / \rho \vert \ge \mu\right\} \right) = 0, \quad \mu >0
$$
that is equivalent to the following 
convergence, 
$$
P\left( \lim_{n\to\infty} {\tilde \rb _1^{(n,\rho)}\over \rho } =0\right)=1. 
$$
Then (3.10) follows. Theorem  3.2 is proved. $\Box$ 
\vskip 0.2cm 
{\it Remark.} Regarding  (3.11) with $\mu = \mu' /\sqrt{\chi_n}$,
we get  the upper bound
$$
P\left( { \vert \tilde \rb_1^{(n,\rho_n)}\vert \over \rho_n} \ge 
{\mu'\over  \sqrt{\chi_n}} \right) 
\le {1\over n^{2(\lfloor \ln n\rfloor / \ln n )\ln(\mu'(1+o(1)))}}, \quad n\to\infty.
$$
Then for any $\mu'>e$ there exists $n_2$ such that 
$$
\sum_{n=n_2}^\infty 
P\left(  \sqrt{\chi_n} 
{\vert \ti \rb^{(n,\rho_n)}\vert \over \rho_n} >\mu'\right)< \infty,
\quad \rho_n= \chi_n \ln n,
$$
and therefore for any 
$\vep >0$ the following relation holds,
 $$
 P\left( \limsup_{n,\chi_n\to\infty}\left\{\omega:  \sqrt { \chi_n} 
 { \vert \tilde \rb_1^{(n,\rho_n)}\vert \over \rho_n} > e + \vep\right\}\right)
 =0, \quad \rho_n= \chi_n \ln n.
 \eqno (3.12)
 $$

\vskip 0.2cm 
Let us consider asymptotic properties of the maximal vertex degree
$d_{\max}^{(n,\rho)}$ (3.1).
\vskip 0.2cm 
{\bf Theorem 3.3.} {\it 
If $\rho_n = \chi \ln n$ with given $\chi >0$, then 
for any  $\ti s > e^{\ti v}$ with $\ti v = \ti h(u)$ given by  
(2.17) the following relation holds,
$$
\lim_{n\to\infty} P\left(  
 \vert {d^{(n,\rho_n)}_{\max}/ \rho_n} - 1\vert \ge   \ti s \right) 
= 0.
\eqno (3.13)
$$
If $\rho_n = \chi_n \ln n$ and $\chi_n \to\infty$ as $n\to\infty$,
then 
$$
P\left( \lim_{n\to\infty}  
 {d^{(n,\rho_n)}_{\max}\over  \rho_n} = 1
 \right) = 1
 \eqno (3.14)
$$
and 
$$
P\left( \limsup_{n \to\infty} \sqrt \chi_n 
{\vert d^{(n,\rho_n)}_{\max} - \rho_n\vert \over \rho_n} > 1\right) = 0. 
\eqno (3.15)
$$

\vskip 0.2cm 
Proof.} Omitting the subscripts in $\rho_n$ and $\chi_n$, we can write that 
$$
\vert \ti d_{\max}^{(n,\rho)} \vert = \vert  d^{(n,\rho)}_{\max} - \rho (n-1)/n\vert \le 
\max_{1\le i\le n} \vert \ti \rb^{(n,\rho)}_i\vert.
 $$
This elementary inequality implies that 
$$
P\left( \vert \ti d_{\max}^{(n,\rho)} \vert > s\right)\le  
P\left(\cup_{i=1}^n \{ \vert \tilde \rb_i^{(n,\rho)} \vert >s\}\right)
\le 
%\sum_{i=1}^n P(\{ \tilde d_i^{(n,\rho)} \vert >s\})= 
n P\left(\vert  \tilde \rb_1^{(n,\rho)} \vert >s\right).
\eqno (3.16)
$$
Regarding inequality (3.3) with $2k $ replaced by $2mk$ and slightly modifying computations of (3.9), we can write that if $\ti s' = e^{\ti v}(1+\delta)$, where $\delta>0$ is such that $\ti s= e^{\ti v} (1+2\delta)$, then (2.20) with $s= \ti s'\rho_n$ implies the upper bound
$$
nP(\vert \tilde \rb^{(n,\rho)}_1\vert \ge  \tilde s' \rho) \le n \left( 
{e^{\ti v}\over \ti s'}(1+o(1)\right)^{2mk}= \exp\left\{ - 2m\lfloor \ln n\rfloor 
\ln(1+\delta) +\ln n \right\}.
$$
There exists $m$ such that 
$2m \ln(1+\delta)  \lfloor \ln n\rfloor /\ln n >1$ for all $n$ starting from certain $n'_0$. 
Then (3.16) implies 
that
$$
\lim_{n\to\infty} P\left( \vert \ti d_{\max}^{(n,\rho)} \vert > \ti s'\rho \right)=0.
$$
Taking into account  obvious estimate
$
\vert {d_{\max}^{(n,\rho)}} - \rho\vert \le \vert 
{\ti d_{\max}^{(n,\rho)}} \vert + 1/n,
$
we get (3.13).

Similarly to (3.11), we can write that for any $\mu>0$
$$
nP\left({\vert \ti \rb_1^{(n,\rho)}\vert \over \rho} \ge \mu\right)
%=P\left( \vert d_1^{(n,\rho)}/ \rho - (n-1)/n \vert \ge \mu\right)
\le {1\over n^{2 \ln(\mu (1+o(1))\sqrt \chi) \lfloor \ln n\rfloor / \ln n -1}}.
$$
Starting from certain $n'_1$, the series of these probabilities converges
because $\chi\to\infty$
and thus (3.14) is true. 
Regarding $\mu = \mu'/\sqrt \chi$ with 
 $\mu'=1+\nu$, $\nu>0$, we get inequality 
$$
n P\left( { \vert \tilde d_1^{(n,\rho)}\vert \over \rho} \ge {\mu'\over \sqrt \chi} \right) \le 
{1\over n^{2m \ln(1+\nu) \lfloor \ln n\rfloor / \ln n -1 }}.
$$
It is clear that  there exists $m$ such that $2m \ln(1+\nu) \lfloor \ln n\rfloor / \ln n>2$ for all $n$
starting from certain $n_1'$. Then (3.15) follows. Theorem 3.3 is proved. $\Box$

\vskip 0.2cm
Let us note that Theorem 3.3 is in deep relation with the Erd\H os-R\'enyi limit theorem \cite{ER},
where the maximum of sums of random variables has been studied,
$$
\Upsilon_n^{(p)} = \max_{i=1,\dots, n} {1\over p} \sum_{j=1}^p w_{ij},
\eqno (3.17)
$$
with $\CW= \{w_{ij}\}_{i,j\in \bN}$ being jointly independent random variables
with zero mean value. It is shown that in the limit $n,p\to\infty$, the value 
$p_c=\ln n$ is critical with respect to the asymptotic behavior of $\Upsilon_n^{(p)}$;
the limiting expression for $\Upsilon_n^{(\kappa \ln n)}$ is obtained \cite{ER}.
In paper \cite{K-02}, similar results  have been obtained in the limit $n,\rho\to\infty$ for random variables 
$$
\CU_n^{(\rho)} = \max_{i=1,\dots, n} U_i^{(n,p)}, 
\quad U_i^{(n,p)}= {1\over \rho} \sum_{j=1}^n a_{ij} w_{ij},
\eqno (3.18)
$$
where $\{a_{ij}\}_{i,j\in \bN}$ is a  family of jointly independent Bernoulli 
random variables with mean $\rho/n$, also independent from  $\CW$.
These random variables are  similar to the normalized maximal vertex degree
(3.2) we have considered,
$$
\CV_n^{(\rho)} = {1\over \rho} d_{\max}^{(n,\rho)} = 
\max_{i=1, \dots, n} {1\over \rho} \sum_{j=1}^n \ra_{ij},
\eqno (3.19)
$$
and Theorem 3.3 
makes the evidence that   the value $\rho = \ln n$ is the critical one with respect to the 
limiting behavior of  $\CV_n^{(\rho)}$. The difference between $\CU_n^{(p)}$ and
$\CV_n^{(\rho)}$ is such that the random variables $\rb_i$ and $\rb_{i'}$ are no more
the independent ones.
One can consider analogous to (3.18) random variables 
$$
Y_n^{(\rho)} = \sum_{j=1}^n a_j^{(n,\rho)} w_j,
\eqno (3.20)
$$
where $a_j^{(n,\rho)}$ are as in (2.1) and $\{w_j\}_{j\in \bN}$ is the 
family of jointly independent identically distributed 
random variables also independent from $\CA^{(n,\rho)}$. In the limit of large $n$, random variables $Y_n^{(\rho)}$ follow the compound Poisson distribution \cite{KPST}. 
In  particular case when $w_j$ take values $\pm 1$ with equal probability,
even moments of $Y_n^{(\rho)}$ are asymptotically  close to the variables (cf. (2.7) and (2.8)),
  $$
\check  {\cal B}_{2k}(\rho) = \sum_{(l_2,l_4, \dots,l_{2k})''}  
\check  B^{(2k)}_{ l_2 ,l_4,\dots, l_{2k}}\,  \rho^{2(l_2+l_4+\dots +l_{2k})},
\eqno (3.21)
$$ 
with 
$$
\check   B^{(2k)}_{l_2,l_4,\dots, l_{2k}} = 
{(2k)!\over (2!)^{l_2} l_2! (4!)^{l_4} l_4! \cdots ((2k)!)^{l_{2k}}l_{2k}!}.
\eqno (3.22)
$$
and the sum in (3.21) 
runs over such integers $l_i\ge 0$
that $l_2+ 2l_4 + \dots +kl_{2k} = k$. 
The value $\check  B_{2k} = \check  \CB_{2k}(1)$ represents the number of partitions
of the set of $2k$ elements into blocks of even cardinality \cite{OEIS}. 
It would be natural to refer to $\{ \check B_{2k}\}_{k\in \bN}$ as to the strongly restricted Bell numbers.
Regarding the asymptotic behavior of strongly restricted Bell polynomials
$\check  {\cal B}_{2k}(\rho)$ in the limit of infinite $k$ and $\rho$,
one can prove analogous to Theorem 2.2 statements. We postpone
the study of the moments of compound Poisson distributions to subsequent publications. 

To complete this section, let us note that the moments $\ti M_k^{(n,\rho)}$ (3.4)
that we study 
can be viewed as a particular case of the moments
$$
L_k^{(n,\rho)}(q) = {1\over n} 
 \sum_{i=1}^n \bE' \left( S_i^{(n,\rho)}(q)\right)^k,
 \quad 
 S_i^{(n,\rho)}(q)= 
 \sum_{j=1}^n\left( \RA_n^{(\rho)}\right)^q_{ij}
 \eqno (3.23)
$$
considered in the case of $q=1$. 
Random variable $S_i^{(n,\rho)}(q)$ counts the number of all possible $q$-step walks
over the graph $\Gamma^{(n,\rho)}$ starting from  vertex $i$. 
It would be interesting to study the asymptotic behavior of 
the moments $L_k^{(n,\rho)}(q)$, $q\ge 2$ (3.23) and the moments of the centered 
random variables $\ti S_k^{(n,\rho)}(q)$ in the limit of infinite
$n,\rho$ and $k$. The limiting values of cumulants of $S_i^{(n,\rho)}(q)$
have been considered in paper \cite{K-08}.

%\newpage
%%%%%%%%%%%%%%%%%%%%%%%%%%%%%%%%%%%%%%%%%%%%%%%%%%%%%%%%%%%%%%%%%%%%%%%%%%%%%%%%%%%%%%
%%%%%%%%%%%%%%%%%%%%%%

\section{Proof of Theorem 2.1 and Theorem 2.2}

To prove Theorem 2.1 and  Theorem 2.2,
we use the method proposed in \cite{TE} to study the asymptotic behavior
of Bell numbers $B_k$. 
It is based on the observation that the local limit theorem holds for 
an auxiliary random variable $Z^{(u)}$ such that $P(Z=k)$ is proportional to 
$u^kB_k/k!$. In the combinatorial studies,  the idea 
to get asymptotic expressions  with the help of the Central Limit Theorem and the local limit theorem
  dates back to the  works by E. A. Bender \cite{Ben1} (see 
 paper \cite{Gawr} and references therein for further developments of the method and 
 monograph \cite{FS} for more detailed information, various applications and 
 generalizations of this approach). 
In these studies of asymptotics  of partitions, main results 
 concern mostly 
 properties of the Stirling numbers of the second kind $S^k_r$, $1\le r\le k$
 as $k,r\to\infty$ and can be applied to the sequence of Bell numbers.
 Further use of this method for the Bell polynomials and restricted Bell polynomials would require proofs
 of a number of additional statements, such as the log concavity of sequences $\CB_k(x)$ 
 and $\tilde \CB_k(x)$ with given $x$ needed in the proof of the local limit theorem for corresponding random variables $Z$. 
 We prefer to stay here within the frameworks of the  approach 
 outlined in \cite{TE}.

In this section, we give the detailed proof of Theorem 2.1
and then describe the modifications of the arguments needed to prove  Theorem 2.2.
Let us first outline the scheme of the proof  based on the method of \cite{TE}. 
We introduce  an auxiliary random variable $Z^{(x,u)}$ that takes values in 
$\bN$,
$$
P(Z^{(x,u)}=k) ={\cal B}_k(x) {u  ^k\over  k! \, G(x,u )}, \quad k\ge 0, \ u>0,
\eqno (4.1)
$$
where 
$
G(x,u ) = \sum_{k=0}^\infty  {\cal B}_k(x) {u ^k/ k!}.
$
Since  $\CB_k(x)$ is equal to the $k$-th moment of the Poisson distribution
$\CP(x)$, we get 
$$
G(x,u ) = \exp\{ x (e^u  -1)\}.
\eqno (4.2)
$$
The generating function  
$F_{x,u } (\tau ) = \sum_{k=0}^\infty  P(Z^{(x,u )}=k) \, \tau ^k$
verifies the following equality
$$
F_{x,u}(\tau) = 
{G(x,\tau u )\over G(x , u)}
\eqno (4.3)
$$
and by elementary computations we conclude that 
$$
{\bf E} Z^{(x,u )} =x ue^u  \quad \mbox{and} \quad   
Var(Z^{(x,u )})= \s_Z^2=  x u(u+1) e^u.
\eqno (4.4)
$$
The local limit theorem that we prove below says that 
if $k= xue^u \to \infty$, then 
$$
P(Z^{(x,u)} = k) = {1\over \sqrt {2\pi} \s_Z}(1+o(1)).
%\eqno (4.4)
$$
This relation, together with  (4.1),
implies asymptotic equality 
$$
\CB_k(x) = {k! \, G(x,u)\over \sqrt {2\pi}\,  \s_Z \, u^k}(1+o(1)), \quad k\to\infty
\eqno (4.5)
$$
and the results of Theorem 2.1 will follow from (4.4) and (4.5).
 
 \subsection{Central and Local Limit Theorems}

Let  
$
\Phi_{Y^{(x,u)}}(t) = \E  e^{-it Y^{(x,u)}},
$
where 
 $$
Y^{(x,u)} =  {Z^{(x,u)}-\E Z^{(x,u)}\over \s_Z^{(x,u)}}, \quad 
\s_Z^{(x,u)} = \sqrt{Var(Z^{(x,u)})}.
$$ 
and $\E Z^{(x,u)}$ and $\sigma_Z^{(x,u)}$ are given by (4.4).  
Given a sequence $(x,u)_N = (x_N,u_N)$, $N\in \bN$, we write that 
 $Y_N= Y^{(x_N,u_N)}$.

 We denote by $({\bf E}Z)_N\to \infty$ the limiting transition such that
 $x_N u_N e^{u_N} \to\infty$ with non-vanishing $x_N$  as $N\to\infty$. 
  Let us consider an  infinite sequences of positive integers 
  $\{k'_N\}_{ N\in \bN}$ and 
  $\{k''_N\}_{N\in \bN}$
 such that 
 $$
  k'_N - x_N u_N e^{u_N} = 
  O\left(\s_Z^{(x,u)}\right), \quad ({\bf E}Z)_N\to\infty
  %\eqno (4.1)
 $$
 and 
 $$
  k''_N - x_N u_N e^{u_N} = 
  o\left(\s_Z^{(x,u)}\right), \quad ({\bf E}Z)_N\to\infty
  %\eqno (4.1)
 $$
   We denote these limiting transitions by $(k',{\bf E}Z)_N\to\infty$
 and $(k'',{\bf E}Z)_N\to\infty$, respectively.
 In what follows, we  omit the superscripts $(x,u)$ 
as well as the subscript $N$ when no confusion can arise.

\vskip 0.2cm
{\bf Lemma 4.1.} 
{\it 
If $(\BE Z)_N\to\infty$, then 
$$
\lim_{N\to\infty} \Phi_{Y}(t) = \exp\{ - t^2/2\}
\eqno (4.6)
$$
for any given $t\in \bR$.
Moreover, 
relation 
$$
P(Z^{(x,u)}=k') -  {1\over \sqrt{ 2\pi } \s_Z }
\exp\left\{ - { \left(k'- \E Z^{(x,u)}\right)^2\over 2 \s_Z^2} \right\} 
= o\left( \sigma_Z^{-1}\right),
\eqno (4.7)
$$
holds in the limit $(k',{\bf E}Z)_N\to\infty$
and 
in particular,  
$$
P(Z^{(x,u)}=k'') =  {1\over \sqrt{ 2\pi } \s_Z }
+ o\left( \sigma_Z^{-1}\right) \quad {\hbox{as}} \quad (k'',{\bf E}Z)_N\to\infty.
\eqno (4.8) 
$$
}

\n{\it Proof.}  Let us first note that
$$
\Phi_Y(t)= e^{-it \E Z/\s_Z}F_{x,u} \left( e^{it/\s_Z}\right).
$$
Relations (4.2) and (4.3) imply equality 
$$
F_{x,u}\left( e^{it /\s_Z}\right) = 
\exp\left\{ xe^u \left(e^{u\Delta} -1 \right)\right\},
\eqno (4.9)
$$
where we denoted $\D = e^{it/\sigma_Z} -1$.
Using expansion
$$
u\D = {iut \over \s_Z} + {u\over 2} \left( {it\over \s_Z}\right)^2 + O\left({ut^3\over \s_Z^3}\right), \quad ({\bf E}Z)_N\to\infty
%\eqno (4.5) 
$$
and observing that
$$
{ut\over \s_Z}= 
{ut\over \sqrt{xu(u+1)e^u}}\le 
{t\over \sqrt {x e^u}} \to 0
\eqno (4.10)
$$ 
in the limit $ ({\bf E}Z)_N\to\infty$, 
we get equality
$$ 
e^{u\Delta} -1 =
 i {ut\over \s_Z} - {t^2\over 2\s_Z^2}u(u+1)
+ O\left({ut^3\over \s_Z^3}\right) + O\left({u^2t^3\over \s_Z^3}
+{u^2t^4\over \s_Z^4} + {u^2 t^5\over \s_Z^5}\right), \quad ({\bf E}Z)_N\to\infty.
$$
Substituting this asymptotic equality   into the right-hand side of (4.9) and remembering
(4.4), we obtain
the following relation,
$$
\Phi_Y(t)= 
\exp\left\{ - {t^2\over 2} + O\left( {xu(u+1) e^ut^3\over \s_Z^3}\right)\right\}
= \exp\left\{ - {t^2\over 2} + O\left( {t^3\over \s_Z}\right)\right\},
\ ({\bf E}Z)_N\to\infty.
\eqno (4.11)
$$
Then (4.6) follows due to (4.10).

\vskip 0.2cm 
We start  the proof of (4.7)  with the computations used   by T. Tao \cite{Tao} in the classical situation of sums of independent random variables.
Taking the mathematical expectation of  the both parts of  identity
$$
{\bf I}_{\{Z = k\}}(\omega) =
{1\over 2\pi} \int_{-\pi}^\pi e^{ipZ} e^{-ipk } dp,
$$
we get by the use of the Fubini's theorem that
$$
P(Z=k) = 
%{1\over 2\pi }\int_{-\pi}^\pi \Phi_{Z}(s) e^{-isk } =
{1\over 2\pi }
\int_{-\pi}^\pi {\bf E}\left( e^{ip(Z-{\E}Z)} \right)e^{-ip(k-{\E}Z)} dp
$$
$$
={1\over 2\pi \s_Z}\int_{-\pi \s_Z}^{\pi \s_Z}
 \Phi_{Y}(t) e^{ -it(k-{\E}Z)/\s_Z } dt,
 \eqno (4.12)
 $$
 where $Y = (Z-{\E Z})/\s_Z$. 
 
 Assuming that the following equality is true
$$
D_N= \int_{-\pi \s_Z}^{\pi \s_Z}
 \Phi_{Y}(t) e^{ -it(k-{\E}Z)/\s_Z } dt
- \int_{-\pi \s_Z}^{\pi \s_Z} 
 e^{ it ( {k-{\E} Z)/ \s_Z} -{t^2/ 2}} dt
 =o(1) 
 \eqno (4.13)
 $$
in the limit $(\E Z)_N\to\infty$, 
we can conclude that (4.7) follows from the identity 
 $$
 {1\over 2\pi} \int_{-\infty}^\infty 
 e^{ i {t(k-{\E} Z)/\s_Z} -{t^2/2}} dt= 
 {1\over \sqrt{2\pi}}e^{ -{(k - \E Z)^2/( 2\s_Z^2)}}
 \eqno (4.14)
 $$
that holds for any sequence $k_N$
and elementary estimate 
$$
 \vert \int_{\vert t \vert >\pi \s_Z} e^{-t^2/2 + i\a t} \vert dt 
 \le \int_{\vert t \vert >\pi \s_Z} e^{-t^2/2}  dt = o(1), \quad \alpha \in \bR.
\eqno (4.15)
$$

To prove (4.13), 
we split the interval of integration into two parts
and
consider first the difference
$$
D_N^{(1)} = 
\int_{\vert t\vert \le y'} \left( \Phi_Y(t) - e^{-t^2/2}\right) e^{-it (k-\E Z)/\s_Z} \ dt,
$$
where $y'$ tends to infinity but  not so fast to guarantee that (4.6) still holds uniformly 
for $t\in [-y',y']$. As we will see, the choice of 
$y' = (\s_Z)^{1/16}$ is convenient for this purpose. 
Using (4.11),  
we conclude  that 
$$
\vert D_N^{(1)} \vert \le 
\int_{\vert t\vert \le y'} e^{-t^2/2}  
\ \vert  \exp\left\{O \left( { t^3\over \s_Z}\right)\right\}
-1\vert 
\ dt = O\left({(y')^4 \over \s_Z}\right) = o(1), \quad (\E Z)_N\to\infty. 
$$

Regarding the difference
$$
D_N^{(2)} = \int_{y'< \vert t\vert \le \pi \s_Z} 
\left( \Phi_Y(t) - e^{-t^2/2}\right) e^{-it (k-\E Z)/\s_Z} \ dt,
$$
we can write that 
$$
\vert D_N^{(2)}\vert \le \int_{y'< \vert t\vert \le \pi \s_Z} \vert 
F_{x,u}\left( e^{it /\s_Z}\right) \vert \ dt + \int_{y'<\vert t\vert \le \pi \s_Z} e^{-t^2/2} dt.
\eqno (4.16)
$$
Elementary computations based on (4.2) and (4.3) show that
$$
\vert F_{x,u}\left( e^{i\b }\right) \vert=
\exp\left\{ xe^u\left( \exp\{ u (\cos \b - 1)\} \cos(\sin(u\beta))- 1
\right)\right\}.
$$
Then clearly,
$$
\vert F_{x,u}\left( e^{i\b }\right) \vert\le \exp\left\{ xe^u\left( \exp\{ u (\cos \b - 1)\} - 1
\right)\right\}.
\eqno (4.17)
$$
The upper estimate of the first integral from the right-hand side of (4.16)
 is based on the fact that 
the argument of the exponent of (4.17) is strictly negative
for all $\beta<0$. Indeed,
using elementary estimate
$$
\cos\left({t\over \s_Z}\right) - 1 
\le - {t^2\over 8 \s_Z^2}, \quad \vert t\vert \le \pi \s_Z,
\eqno (4.18)
$$
and its consequence 
$$
\cos\left( {t\over \s_Z} \right) - 1\le - {(y')^2\over 8\s_Z^2},
\quad y'\le \vert t\vert\le \pi \s_Z,
$$
we conclude that 
$$
\exp\left\{ u \left(\cos \left({t\over  \s_Z}\right) -1\right)\right\} \le 
\exp\left\{ - {u(y')^2\over 8\s_Z^2}\right\}.
%\eqno (4.18)
$$
Observing that in the limit $(\E Z)_N\to\infty$, the argument of the last exponential
tends to zero, we apply to it the Taylor expansion and write that for sufficiently large $xue^u$, 
$$
x e^u \left( \exp\left\{ - {u(y')^2\over 8\s_Z^2}\right\}
- 1\right) 
\le -x u e^u { (y')^2\over 8\s_Z^2} = - {(y')^2\over 8(u+1)}
$$
as $(\E Z)_N\to\infty$. Then 
we deduce from (4.17) that 
$$
\vert F_{x,u}\left( e^{it/\s_Z }\right) \vert\le \exp\left\{- {(y')^2\over 8(u+1)}
\right\}, \quad (\E Z)_N\to\infty
$$
and therefore
$$
\int_{y'< \vert t\vert \le \pi \s_Z} \vert 
F_{x,u}\left( e^{it /\s_Z}\right) \vert \ dt
\le 2\pi \s_Z \exp\left\{- {(y')^2\over 8(u+1)}\right\}
=2\pi \s_Z \exp\left\{- {\s_Z^{1/8}\over 8(u+1)}\right\}
.
\eqno (4.19)
$$
Equality (2.12) shows that in the asymptotic regimes when either 
$x=O(k)$ or $x/k\to\infty$ as $k\to\infty$, the parameter $u$
either is finite or tens to zero. Then
it is clear that the right-hand side of (4.19) vanishes 
in the limit \mbox{$(\E Z)_N\to\infty$}  and (4.16) implies relation
$\vert D_N^{(2)}\vert = o(1)$ because $y' = \s_Z^{1/16} \to\infty$.
Then  (4.13) follows. 

If $x= o(k)$, then (2.12) implies relation
$u\approx \ln(k/x)$ and $\s_z \approx xu^2 e^u= uk$. Then 
elementary analysis shows that the right-hand side of (4.19)
vanishes, with the choice of $y' = \s_Z^{16}$, 
under additional condition that  $k^{1/7} \gg \ln(k/x_k)$
(see Remark after Theorem 2.1). 
In Theorem 2.1, as well as in the definition of 
limiting transition $(\E Z)_N\to\infty$, 
we assume that the sequence $(x_k)_{k\in \bN}$ 
does not converge to zero and therefore this additional condition is obviously verified and the right-hand side of  (4.19) goes to zero as 
$k\to\infty$.

\vskip 0.2cm 
Gathering (4.12), (4.13), (4.14) and (4.15),
we conclude that (4.7) is true and, 
 in particular, (4.8) holds. Lemma 4.1 is proved. $\Box$

%%%%%%%%%%%%%%%%%%%%%%%%%%%%%%%%%%%%%%%%%%%%%%

\subsection{Proof of Theorem 2.1}

Equation 
$
ue^u = \beta, \ \beta>0
$
(2.12)
has a unique solution $u= u(\beta)$ known as the Lambert $W$ function 
\cite{DB,FS}. Given an infinite sequence $\left\{ (x_k)_{k\in \bN}\right\}$ of strictly positive reals, 
we determine 
$u_k$ such that 
$$
u_k e^{u_k} = {k\over x_k}, \quad k\in \bN.
%\eqno (4.19)
$$
In this subsection,  we   omit the subscripts $k$ in $x_k$ and $u_k$
when no confusion can arise.

Rewriting (4.1) in the form  
$$
{\cal B}_k(x) = P(Z^{(x,u)} = k) {k!\over u^k} G(x,u),
$$
we get with the help of (4.8) the following asymptotic equality,
$$
{\cal B}_k(x) = {1\over  \sqrt{ 2\pi x u(u+1) e^u} }   \exp\{ x (e^u-1)\} \, {k!\over u^k}\, (1+o(1)), \quad k\to\infty.
%\eqno (2.16)
$$ 

Using the Stirling formula (2.22)
and (4.19), we get  relation   
$$
{\cal B}_k(x)= {x^k\over \sqrt{u+1}} \exp\left\{ x u(u-1) e^u + x(e^u-1)\right\} (1+o(1)), \quad k\to\infty,
$$
and finally 
$$
 {\cal B}_k(x) = {x^k\over \sqrt {u+1}} 
\exp\left\{ k \left( u -1 +{1\over u}\right) - x\right\} (1+o(1)), \quad 
k\to\infty.
\eqno (4.20)
$$

Asymptotic equality (4.20) coincides with the result  
by D. Dominici \cite{D} obtained with the help of the ray method
applied to the differential-difference equation 
$$
\CB_{k+1} (x) = x\left( \CB_k'(x) + \CB_k(x)\right)
\eqno (4.21)
$$
(see also \cite{E}).
Relation (4.20)  considered at $x=1$ gives an expression 
for the Bell numbers $B_k= \CB_k(1)$ 
 similar to that obtained by E. G. Tsylova and E. Ya. Ekgauz \cite{TE}.

Now we will examine the asymptotic behavior of the sequence  ${\cal B}_k(x_k)$, $k\to\infty$  in dependence  whether  $0< x_k\ll k$, or $x_k=O(k)$, or  
$x_k\gg k$.
Regarding an  auxiliary variable,  
$
\Psi_k(x)= k^{-1}  \ln \left( {{\cal B}_k(x)/ x^k} \right),
$
we deduce  from (4.20) that 
$$
\Psi_k(x) = u - 1 + {1\over u} - {1\over u e^u} - {1\over 2k} \ln (u+1) + o(k^{-1}), \quad k\to\infty.
\eqno (4.22)
$$
\vskip 0.2cm 
\noindent 
{a)} If $x_k/k\to 0$, then the right-hand side of  (2.12) 
tends to infinity. It is not hard to see that   the solution $u= u(\beta)$ of the transcendent equation 
$ue^u=\beta$ has the following asymptotic expansion
\cite{DB},
$$
u = \ln \beta  - \ln \ln {\beta}   + 
O\left( { \ln \ln \beta  \over \ln \beta }\right), 
\quad \beta \to \infty.
\eqno (4.23)
$$
Substituting this expression with $\beta = k/x$  into the right-hand side 
of (4.22), we get
the following  
asymptotic equality, 
$$
\Psi_k(x) = \ln \left({k\over x}\right) - \ln \ln \left({k\over x}\right) - 1 + 
O\left({\ln \ln \left( {k/ x}\right)\over \ln(k/x)}\right),
 \quad {k\over x} \to \infty, 
\quad k\to\infty.
$$
Returning to the variable ${\cal B}_k(x)$,
we can write that if $x/k\to 0$, then 
$$
{\cal B}_k(x) = x^k \exp\left\{ k \ln \left( { k\over x}\right)  - k \ln \ln \left( { k\over x}\right) - k  + O\left({k \ln \ln \left( {k/ x}\right)\over \ln(k/x)}\right)
\right\}, \ k\to\infty.
\eqno (4.24)
$$
This relation implies (2.9). 

Regarding (4.23)  in the particular case   $x=1$, we get the following asymptotic equality 
for the Bell numbers $B_k = \CB_k(1)$,
$$
{\ln B_k\over k} = \ln k - \ln \ln k - 1 + O\left( { \ln \ln k \over \ln k}\right),\quad 
k\to\infty
%\eqno (2.23)
$$
that is equivalent to the result of  \cite{TE}. 
The first three terms of the right-hand side of this relation 
reproduce  those of the  asymptotic expansion
of the  Bell numbers obtained by N. G. de Bruijn \cite{DB}
(see also  papers \cite{Lov} and \cite{MW}).

\vskip 0.2cm 
\noindent {b)} If $\chi = x_k/k$ as $k\to\infty$,  
then  relation (4.22) implies  equality
$$
\lim_{k\to\infty} \Psi_k(x_k) = u- 1 +{1\over  u} - {1\over u e^u} = h(u), \quad ue^u= {1\over \chi}.
\eqno (4.25)
$$
Relation (2.10) follows from (4.25) with $v=h(u)$ given by (2.11) (see also (2.23)).
If $(x_k)_{k\in \bN}$ is such that $x_k = \chi k$ for all $k\in \bN$, then (2.13)
follows directly from (4.20).
\vskip 0.2cm

\n {c)}  Consider the last asymptotic regime when 
$ x_k/k = \chi_k \to\infty $, $k\to\infty$.  
In this case  $u\to 0$ (4.19) 
and 
$$
u = {1\over \chi_k} - {4\over \chi^2_k} + o(\chi^{-2}_k), \quad k\to\infty.
$$
In this case 
$
\CH_k(x) = {1\over u} \left(1 - {1\over e^u}\right) + u - 1 + o(1/k) = {u\over 2} + o(u^2) + o(1/k), \ k\to\infty.
$
Then 
$$
\CH_k(x) = {1\over 2\chi_k} - {2\over \chi^2_k}  + o(\chi^{-2}_k) + o(k^{-1}), \quad k\to\infty 
$$ 
and therefore
$$
{\cal B}_k(x) = \left(k \chi_k  \exp\left\{{1\over 2\chi_k} - {2\over \chi^2_k} + o(1/k)\right\}\right)^k, 
\quad   k\to\infty.
%\eqno(4.26)
$$
Then (2.14) follows. Theorem 2.1 is proved. $\Box$

%%%%%%%%%%%%%%%%%%%%%%%%%%%%%%%%%%%%%%%%%%%%%%%%%%%%%%%%%
\subsection{Restricted Bell polynomials and proof of Theorem 2.2}

In Section 5, we show that the  exponential generating function
of the sequence $\{\tilde \CB_k(x)\}_{k\ge 0}$ is given by
$$
\tilde G(x,u) 
= \exp\{ x(e^u - u-1)\}. 
\eqno (4.26)
$$
Regarding  a random variable $\tilde Z^{(x,u)}$ such that   
$$
P( \ti Z^{(x,u)} = k) =  \ti \CB_k(x) {u^k\over k!\,  \tilde G(x,u)},
\quad k\ge 0,
$$
it is easy to see that 
the generating function $\tilde F_{x,u}(\tau)= \sum_{k\ge 0} \tau^k \ti p_k$
is given by equality 
$$
\tilde F_{x,u}(\tau) = { \tilde G(x,\tau u)\over \tilde G(x,u)}.
\eqno (4.27)
$$ 
Elementary computations based on (4.27) show that 
$$
\E \ti Z^{(x,u)} = xu(e^u-1)\quad \mbox{and} \quad 
Var(\ti Z^{(x,u)}) = \ti \s^2= xu((u+1)e^u -1).
$$
We introduce a random variable
$
\ti Y^{(x,u)} = { \ti Z^{(x,u)}-\bE \ti Z^{(x,u)}/ \ti \s}
$
with the characteristic function
  $$
  \Phi_{\ti Y} (t) = \E e^{it \ti Y} = e^{-i t \bE \ti Z/\ti \s} 
  \ti F(e^{i t /\ti \s}).
  \eqno (4.28)
  $$ 
  Here and below, we omit the superscripts $(x,u)$ in $\ti Y^{(x,u)}$ and
  $\ti Z^{(x,u)}$ when no confusion can arise. 
  
 Let us  consider a sequence $(x_N,u_N)_{N\in \bN}$ such that 
 $(x_N)_{N\in \bN}$ does not converge to zero,
 $$
 x_N u_N\left( e^{u_N} -1\right) \to\infty
 \quad {\hbox{and}} \quad  x_N (u_N\left( e^{u_N} -1\right))^3 \to\infty
 \quad {{\hbox{as}}}\ \ N\to\infty
\eqno (4.29)
$$
 and denote this limiting transition  by $(\E\tilde Z)^{(2)}_N\to\infty$.

  \vskip 0.2cm 
\noindent   {\bf Lemma 4.2.} {\it  
  If $\ti k_N$ is such that 
 $$
 \ti k_N - x_N u_N\left( e^{u_N} -1\right) = o(\ti \s_N), \quad (\E\tilde Z)^{(2)}_N\to\infty,
 $$
 where $\ti \s_N^2 = x_Nu_N((u_N+1)e^{u_N} -1)$ verifying (4.29), then 
 $$
 P\left( \ti Z_N = \ti k_N\right)=
 {1\over \sqrt{2\pi} \ti \s_N} (1+o(1)),
 \eqno (4.30)
 $$
 where we denoted $\ti Z_N= Z^{(x_N,u_N)}$.

 \vskip 0.2cm
 Proof.
 } Let us first  show that  the Central Limit Theorem holds for
 the random variables $\ti Y^{(x,u)}$,  
 $$
 \ti \Phi_N(t)= 
 \Phi_{\ti Y_N} (y) =  e^{-t^2/2}(1+o(1)), \quad (\E \ti Z)^{(2)}_N\to\infty,
  \eqno (4.31)
  $$
 where $\ti Y_N = \ti Y^{(x_N,u_N)}$. In what follows, we omit the subscripts
 $N$ when no confusion can arise.
It follows from (4.27) and (4.28)  
that 
$$
\ti F\left(e^{it/\ti \s }\right) = 
\exp\left\{ x e^u \left( e^{u\ti \Delta} - 1\right)- xu\ti \Delta
\right\}, \quad \tilde \Delta = e^{it/\ti \s} - 1.
$$
Similarly to (4.11), we can write that
$$
u\ti \Delta = {iut\over \ti \s} - {ut^2\over 2\ti \s^2} + 
O\left({ut^3\over \ti \s^3}\right)
$$ 
and 
$$
e^{u\ti\Delta} -1 =
 i {ut\over \ti\s} - {t^2\over 2\ti \s^2}u(u+1)
+ O\left({ut^3\over \ti\s^3}\right) + O\left({u^2t^3\over \ti\s^3}+
{u^2t^4\over \ti\s^4} 
+ {u^2 t^5\over \ti\s^5}\right), \quad ({\bf E}\ti Z)^{(2)}_N\to\infty.
$$
Then 
$$
\ln \ti F(e^{iy /\ti \s}) = 
{xu(e^u-1) it\over \ti \s} -{t^2\over 2} + 
O\left({xu(u+1)e^ut^3\over \ti \s^3} \right), \quad ({\bf E}\ti Z)^{(2)}_N\to\infty
\eqno (4.32)
$$
Taking into account inequality $\tilde \s^2 \ge  xu^2 e^u$, we conclude that 
$$
{xu(u+1)e^ut^3\over \ti \s^3} \le  {(u+1)t^3\over \sqrt{ x u^4 e^u}}.  
$$
It is easy to see that the right-hand side of this inequality vanishes
in the cases when either $u=O(1)$ or $u\to\infty$ as $xu(e^u-1)\to\infty$. 
If $u\to 0$, then (2.18) implies asymptotic equivalence $u^2\approx k/x$
and the right-hand side of this inequality is of the order 
$t^3/ \sqrt {k^2/x}$. It follows from the second condition of (4.29) 
that in this asymptotic regime $x= o(k^{3/2})$ and then (4.31) follows.

To prove (4.30), we repeat the proof of Lemma 4.1 and see that 
it suffices to show that (cf. (4.13))
$$
\ti D_N= \int_{-\pi \ti \s } ^{\pi \ti \s}
\vert \ti\Phi_N(t) - e^{-t^2/2}\vert dt = o(1), \quad 
({\bf E}\ti 
Z)^{(2)}_N\to\infty.
$$
It follows from (4.28) and (4.32) that 
$$
\vert \ti D_N^{(1)}\vert  \le \int_{\vert t\vert \le \ti y' } e^{-t^2/2} \
\vert \exp\left\{O\left( {xu(u+1)e^ut^3\over \ti \s^3} \right)\right\} -1
\vert \ dt = O\left( {\ti y^4\over \sqrt{xu^2 e^u}}\right) 
+ O\left( {\ti y^4\over \sqrt{xu^4 e^u}}\right). 
$$
Choosing 
$$
\ti y = \min\left\{ (xu^2e^u)^{1/16}, (xu^4e^u)^{1/16} \right\},
\eqno (4.33)
$$ 
we make $\tilde D_N^{(1)}$ vanishing in the limit 
$({\bf E}\ti 
Z)^{(2)}_N\to\infty$.

To estimate the remaining part of $\ti D_N$, we consider 
$\ti F(e^{i\beta})$ and write that 
$$
\vert \ti F (e^{i\beta})\vert  = \exp\left\{ x  
e^u\left( e^{u(\cos \b - 1)}\cos(u\sin \beta)\right)\right\} \exp\{ -xe^u\} 
\exp\left\{ - x u(\cos \b -1)\right\}.
$$
Then obviously
$$
\vert \ti F (e^{i\beta})\vert\le 
\exp\left\{ x  
e^u\left( e^{u(\cos \b - 1)}-1\right)\right\}.
%\eqno (4.33)
$$
Using (4.18), we conclude that  the following upper bound
$$
\vert \ti F (e^{it/\tilde \s })\vert\le
\exp\left\{ - xu e^u { \tilde y^2 \over 16 \tilde \s^2}\right\}
= \exp\left\{ - {e^u \tilde y^2\over 16 ( u e^u + e^u - 1)}\right\}
,
\quad \tilde y \le \vert t\vert \le \pi \tilde \s
%\eqno (4.33) 
$$
is true for sufficiently large values of  $xu (e^u - 1)$. 
Then, in complete analogy with (4.16), 
$$
\vert \tilde D_N^{(2)}\vert \le 
2\pi \tilde \s \exp\left\{ - {e^u \tilde y^2\over 16 ( u e^u + e^u - 1)}\right\} + \int_{\vert t\vert \ge \tilde y} e^{-t^2/2} dt.
\eqno (4.34)
$$
Elementary analysis shows that the first term of the right-hand side of (4.34) vanishes in the case when $u$ remains non-zero and finite as 
$k\to\infty$. If $u\to 0$ as $k\to\infty$, then the first terms of the right-hand side of (4.34) vanishes due to the second condition of (4.29). 
Finally, if $u\to\infty$ and the $\liminf_{k\to\infty} x_k>0$,
then the first term of the right-hand side of (4.34) vanishes. 
This observation completes the proof of Lemma 4.2. $\Box$

Let us note that if $x_k\to 0$ as $k\to\infty$, then 
it is not hard to show that the first term of the right-hand side of (4.34)
vanishes in the limit  
$({\bf E}\ti 
Z)^{(2)}_N\to\infty$
provided $x\gg k e^{- k^{1/15}}$ as $k\to\infty$ (see Remark after Theorem 2.2). We do not present 
the details of this reasoning here also because this asymptotic regime 
does not belong to the main subject of the present paper.

\vskip 0.2cm 

Let us complete the proof of Theorem 2.2.
We consider an infinite sequence 
of strictly  positive reals $(x_k)_{k\in \bN}$ and determine 
$(u_k)_{k\in \bN}$ 
such that $ u_k (e^k-1) = k/x_k$ for all $k\ge 1$.
For these values of $x_k$ and $u_k$,
relation (4.30) implies  that  
$$
\ti \CB_k(x) = {1\over \sqrt{ 2\pi xu((u+1) e^u-1)}} \sqrt{2\pi k} 
\left( {k\over e u }\right)^k e^{x(e^u -u-1)} (1+o(1)),
\quad k\to\infty.
\eqno (4.35)
$$
Here and below we omit the subscripts $k$ in $x_k$ and $u_k$ when no confusion can arise.

Analysis of (4.35) is similar to that performed in the proof of Theorem 2.1. 
In particular, to get relation (2.15), we use an observation that  the 
asymptotic expansion  of the solution $\tilde u= \tilde u(\beta)$
of equation 
$$
\tilde u\left(e^{\tilde u}-1\right) = \beta, \quad \beta\to\infty
\eqno (4.36)
$$
coincides with the right-hand side of (4.23)
  (see Section 5, Lemma 5.4). Then (2.15) follows from relation (4.24).
  Asymptotic equality (2.16) is a direct consequence of (4.35)
and expression for $\ti h(u)$ (2.17) (see also (2.24)). 
Relation (2.19) is a direct consequence of (4.35) considered with  
$x_k = \chi k$.

Finally, to prove   
(2.20), we observe that in this asymptotic regime   
$$
\tilde \CH_k(x_k) = {1\over k} \ln \left(
{\tilde \CB_k(x_k)\over x^k}\right) = \ln(e^u - 1) - 1 +{x_k\over k} \left( e^u - u - 1\right) +o(k^{-1}),
$$
where $u= u_k$ verifies equation
$$
u(e^u-1 ) = {1\over \chi_k}, \quad {1\over \chi_k} =  {k\over x_k } \to 0
$$
and therefore   $u= \sqrt{k/x_k} (1+o(1))$, $k\to\infty$.
Elementary computations show that 
$$
\ln (e^u-1) -1 = \ln \left( u+ {u^2\over 2} +o(u^2)\right) - 1 =
\ln u - 1 - u/2 + o(u)
$$ 
and 
$
x \left( e^u - u - 1\right)/k = {1/2} +o(1)$, as $ k\to\infty$.
Then 
 $$
 \tilde \CH_k(x) = {1\over 2} \left( \ln\left( {k\over x}\right) - 1\right) + o(1),
 \quad   k\to\infty,
 \eqno (4.37)
 $$
and 
$$
\tilde \CB_k(x) = \left( x \exp\left\{ { \ln(k/x) -1\over 2} +o(1)\right\}\right)^k, \quad k\to\infty.
\eqno (4.38)
$$
This gives  (2.20). Additional restriction $x= o(k^2), k\to\infty$
is a consequence of the condition that $(\E \ti Z)^2/x_N \to \infty$ as $N\to\infty$  
imposed in  Lemma 4.2 (see (4.29)).

Theorem 2.2 is proved. $\Box$ 
%%%%%%%%%%%%%%%%%%%%%%%%%%%%%%

%\subsection{Deviation probability of vertex degree }

%%%%%%%%%%%%%%%%%%%%%%%%%%%%%%%%%%%%%%%%%%%%%%%%%%%%
%%%%%%%%%%%%%%%%%%%%%%%%%%%%%%%%%%%%%%%%%%%%%%%%%%%%%
\section{Auxiliary facts and discussion}

In this section we collect  proofs of the statements we have used above 
and formulate  a number of important supplementary observations. 

\subsection{Binomial and Poisson random variables}

Let us describe convergence of
random variables $X_n^{(\rho)}$ in the limit $n,\rho\to\infty$.
We denote by
$
\Phi_{n,\rho}(t) = {\bE} \exp\{ it X_n^{(\rho)}\}
$
the characteristic function of $X_n^{(\rho)}$.
\vskip 0.2cm {\bf Lemma 5.1.}   
{\it If $\rho= o(\sqrt n)$ when $n$ infinitely increases, then 
$
\Phi_{n,\rho}(t)$ converges
to the one of the Poisson distribution, 
$
\Phi_{Y^{(\rho)}}(t) = \exp\{ \rho (e^{it}-1)\}
$
in the sense that for any $t\in \bR$
$$
\Phi_{n,\rho}(t)/\Phi_{Y_\rho} (t) \to   1, \quad n\to\infty.
\eqno (5.1) 
$$
If $k= o(\sqrt n)$ and $\rho = o(\sqrt n)$ when $n$ infinitely increases, 
then 
$$
P(X_n^{(\rho)} = k ) / P(Y_\rho = k) \to  1, \quad n\to\infty.
\eqno (5.2)
$$
}
Some of these results are  known but we formulate them for completeness. 

{\bf Lemma 5.2.}
{\it If $\rho = o(\sqrt n)$ when $n$ tends to infinity, then 
for any given $t\in \bR$
$$
\lim_{n,\rho\to\infty}\E \exp\{i t U_n^{(\rho)}\} = e^{-t^2/2}, 
\eqno (5.3)
$$ 
where 
$
U_n^{(\rho)} = {(X_n^{(\rho)} - \rho)/ \sqrt \rho}
$.
}

The proofs of relations (5.1) and  (5.2) 
 are based on simple use  the Taylor expansions
of characteristic functions. 
Indeed, assuming $\rho= o(n), n\to\infty$, we can write that 
$$
{\bE} \exp\{ it X_n^{(\rho)}\} = \left( e^{it} {\rho\over n} + 
\left(1-{\rho\over n}\right)\right)^n
$$
$$
=\exp\left\{ n\ln \left( 1+ { (e^{it}-1)\rho\over n}\right) \right\}
= \exp\left\{ (e^{it}-1)\rho + O(\rho^2/n)\right\}.
$$
Then (5.1) follows.

Regarding the probability distribution of  $X_n^{(\rho)}$, we can write that 
$$
P(X_n^{(\rho)}=k)= R(k,n)\, 
{\exp\left\{n \ln (1-\rho/n)\right\}\over \exp\left\{ k \ln (1-\rho/n)\right\}} \cdot  {\rho^k\over k!} ,
\eqno (5.4)
$$
where we denoted
$$R(k,n) =  \prod_{i=1}^{k-1} {n-i\over n} = \prod_{i=1}^{k-1} \left( 1 - 
{i\over n}\right).
$$
If $k/n\to 0$, then 
$$
\ln R(k,n) = \sum_{i=1}^{k-1} \ \ln \left( 1-{i\over n}\right) 
$$
$$= 
\sum_{i=1}^{k-1} \left( -{i\over n} + {i^2\over 2n^2} + 
O\left({i^3\over n^3}\right) \right) =- {(k-1)k\over 2n} + 
O\left({k^3\over n^2}\right).
\eqno (5.5)
$$
Using the Taylor expansion of $\ln(1- \rho/n)$, one can easily deduce from 
(5.4) with the help of (5.5) the following relation, 
$$
P(X_n^{(\rho)}=k) {e^\rho k!\over \rho^k} =
\left( 1+ O\left({k^2\over n}\right)\right) 
\left( 1+ O\left({\rho^2\over n}\right)\right) 
\left( 1+ O\left({k\rho\over n}\right)\right). 
$$
This  (5.2) follows.

The proof of Lemma 3.2 is 
elementary and we do not present it here.

\subsection{Generating function of restricted Bell polynomials }

In this subsection we prove the following statement.

\vskip 0.2cm 
 {\bf Lemma 5.3.} {\it The exponential generating function of  the sequence 
 $\{ \tilde \CB_k(x)\}_{k\in \bN}$
 (2.7),  (2.8) is given by relation
 $$
 \tilde G(x,u) = \sum_{k=0}^\infty 
 \tilde \CB_k(x) {u^k\over k!} = \exp\left\{ x\left( e^u - u - 1\right)\right\}. 
 \eqno (5.6)
 $$

 Proof.}
 Let us consider  analogs of the Stirling numbers of the second kind,
$$
\tilde S^k_r= {1\over r!} \sum_{(h_1,h_2, \dots, h_r)'}
{k!\over h_1!\,  h_2! \cdots h_r!}=
$$
$$
= {1\over r!} \sum_{(h_1,h_2, \dots, h_r)'}
 {k \choose h_1} {k-h_1\choose h_2} \cdots {k-h_1-h_2\dots - h_{r-1}\choose h_r},
 \eqno (5.7)
$$ 
 where the sum over $(h_1,h_2, \dots, h_r)'$ is such that $h_1+\dots +h_r = k$ and $h_i\ge 2, i=1, \dots, r$. 
 It is easy to deduce from this definition  that
 $$
 \sum_{k=r}^\infty \tilde S^k_r {t^k\over k!} = {1\over r!} 
 \left( e^t-t-1\right)^r.
 $$
 Taking into account definition  
 $
 \tilde {\cal B}_k(x) = \sum_{r=0}^k \tilde S^k_r x^r,
 $
 we conclude that 
 $$
\ti G(x,t) =  \sum_{k=0}^\infty \tilde {\cal B}_k(x) {t^k\over k!} = 
 %\sum_{k=0}^\infty \sum_{r=0}^k \tilde S(k,r) x^r {t^k \over k!}
 %=
\sum_{r=0}^\infty \sum_{k=r}^\infty \tilde S^k_r x^r {t^k \over k!} 
 = \exp\left\{ x\left(e^t - t- 1\right)\right\},
 $$
 where we have interchanged the order of summation by  standard arguments.
 The last equality completes the proof of (5.6). $\Box$
 
 It follows from (5.7) that  restricted Stirling numbers of the second kind
 verify  recurrence
 $$
 \ti S^{k+1}_r = r \ti S^k_r + k \ti S^{k-1}_{r-1}, \quad 1\le r\le k
 \eqno (5.8)
 $$
 with  obvious initial conditions
 $\ti S^k_0= \delta_{k,0}$, $\ti S^1_1=0$ and $\ti S^k_{k-1+l}=0$, $l\ge 0$.
 It would be interesting to study asymptotic properties of $\tilde S^k_r$
 in the limit $k,r\to\infty$, but this questions is out of  the frameworks of the present paper. 
 
 %%%%%%%%%%%%%%%%%%%%%%%%%%%%
 
 \subsection{Proof of Lemma 2.1 and Lemma 3.1}
 
 %Moments $M_k^{(n,\rho)}$ and $ \tilde  M_k^{(n,\rho)}$ %

To study the moments $\CM_k^{(n,\rho)} =\E\left( \sum_{j=1}^n a_j\right)^k$ 
(2.2), 
it is natural to represent  the multiple sum of the right-hand side of 
this equality 
as a sum over  classes of equivalence $\cal C$, 
each class being associated with a partition of the set $\{j_1,j_2,\dots, j_k\}$ into  blocks
such that  the variables in each block
are equal to  the same value from the set $\{1,2,\dots,n\}$. 
It is easy to see that 
$$
{\bE}\left(  \sum_{j_1=1}^n   \cdots \sum_{j_k=1}^n a_{j_1} a_{j_2}\cdots a_{j_k} \right) 
= \sum_{\{{\cal C}\}} \prod_{i=1}^k
\left( \E a_1^i\right)^{l_i} n(n-1)\cdots (n-\vert {\cal C}\vert  +1),
\eqno (5.9)
$$
where $\vert {\cal C} \vert= l_1+l_2+\dots+l_k$ denotes the number of blocks 
in the partition $\cal C$. 
Here and below we omit the superscripts $(n,\rho)$.
Since $\bE a_1^i= \rho/n$ and $\vert{ \cal C}\vert\le k$, then 
 the elementary estimate (cf. (5.5))
$$
\log \prod_{i=1}^{\vert \CC\vert-1} \left( 1 - {i\over n}\right) = 
- {\vert \CC\vert (\vert \CC\vert -1)\over 2n} + O(\vert \CC\vert^3/n^2), 
\quad n\to\infty
\eqno (5.10)
$$
shows that   if $k=o( \sqrt n)$, then 
$$
\sum_{\{{\cal C}\}} \left( {\rho\over n} \right)^{\vert {\cal C} \vert } n(n-1)\cdots (n-\vert {\cal C}\vert  +1)
= \sum_{\{{\cal C}\}}\ \rho^{\vert {\cal C} \vert}\, (1+o(1)),\quad 
n\to\infty
$$
Therefore in this limit,
$$
\CM_k^{(n,\rho)} =  \sum_{\{{\cal C}\}}\ \rho^{\vert {\cal C} \vert}\, (1+o(1)), \quad n\to\infty
$$
and relation (2.3) follows from the fact that the number of classes $\CC$ with given 
$(l_1,l_2,\dots, l_k)$
is equal to the number $B^{(k)}_{l_1,l_2, \dots, l_k}$ (2.6).

Let us consider the moments $\tilde  \CM_k^{(n,\rho)}$ (2.2). As in (5.9), we have 
$$
 \tilde \CM_k^{(n,\rho)}=
  %{\bf E}\left( \sum_{j_1,j_2,  \dots, j_k=1}^n  \tilde a_{j_1}  \tilde  a_{j_2}\cdots  \tilde a_{j_k} \right) = 
  \sum_{\{{\cal C}^*\}} \prod_{i=2}^k
\left( \E  \tilde a_1^i\right)^{l_i} n(n-1)\cdots (n-\vert {\cal C^*}\vert  +1),
\eqno  (5.11)
$$
where the sum runs over the classes of equivalence $\cal C^*$ given by  such partitions of the set $\{1,2,\dots , n\}$
that have no blocks of one element. 
It is easy to see that
$$
\E \ti a_1^m= \E \left(a_1 - {\rho\over n}\right)^m = {\rho\over n} Q_m(\rho/n),
$$
where 
$$
Q_m(\rho/n)= \sum_{l=2}^m {m \choose l} \left( - {\rho\over n}\right)^{m-l} + (m-1) \left( - {\rho\over n} \right)^{m-1}.
$$
Then obviously
$$
Q_m(\rho/n) = \left( 1-{\rho\over n}\right)^m + \left( {-\rho\over n}\right)^{m-1} \left( 1+ {\rho\over n} \right) 
\le \left( 1 + {2\rho\over n} \right)^m.
\eqno (5.12)
$$
Substituting upper bound (5.12) into (5.11), we can write  that 
$$
\tilde \CM_k^{(n,\rho)}\le  \sum_{\{{\cal C}^*\}} \prod_{i=2}^k
\left({\rho\over n}  \left( 1 + {2\rho\over n}\right)^i \right)^{l_i} n(n-1)\cdots (n-\vert {\cal C^*}\vert  +1)
$$
$$
\le \left( 1+{2\rho\over n}\right)^k  \sum_{\{{\cal C}^*\}} 
 \rho^{\vert \CC^*\vert }\, { (n-1)(n-2) \cdots (n- \vert {\cal C}^*\vert +1)\over n^{\vert {\cal C}^*\vert -1}}
.
\eqno (5.13)
$$
Using (5.10), we conclude that if $k= o(\sqrt n)$ and $\rho = o(n)$
as $n\to\infty$, then 
$$
\tilde  \CM_k^{(n,\rho)}\le  \sum_{ \{\CC^*\}} \rho ^{\vert \CC^*\vert} (1+o(1))= 
 \tilde  
 {\cal B}_k(\rho) (1+o(1)), \quad n\to\infty.
 \eqno (5.14)
$$
Elementary analysis shows that the lower estimate  
$$
Q_m(\rho/n) \ge \left( 1-{\rho\over n}\right)^m \left( 1- {4\rho\over n-\rho}\right)  \ge 
 \left( 1-{\rho\over n}\right)^m 
\left( 1- {4\rho\over n-\rho}\right) ^m
$$
is true for $m\ge 2$ and sufficiently large $n,\rho$ such that $\rho = o(n)$. 
Then (5.11) implies inequality
$$
\tilde  \CM_k^{(n,\rho)} \ge  \left( 1-{\rho\over n}\right)^k 
\left( 1- {4\rho\over n-\rho}\right)^k \tilde \CB_k(\rho)
= \tilde \CB_k(\rho)(1+O(\rho k/n)).
\eqno  (5.15)
$$
Relations (5.14) and (5.15) prove the second part of 
Lemma 2.1 given by (2.4).

\vskip 0.2cm 
Let us comment on the proof of Lemma 3.1. 
Similarly to (5.11), we have relation 
$$
 \tilde M_k^{(n,\rho)}=
  %{\bf E}\left( \sum_{j_1,j_2,  \dots, j_k=1}^n  \tilde a_{j_1}  \tilde  a_{j_2}\cdots  \tilde a_{j_k} \right) = 
  \sum_{\{{\cal C}^*\}} \prod_{i=2}^k
\left( \E  \tilde a_1^i\right)^{l_i} (n-1)(n-2)\cdots (n-\vert {\cal C^*}\vert )
$$
that implies the following bounds,
$$
 \tilde M_k^{(n,\rho)}
\le \left( 1+{2\rho\over n}\right)^k  \sum_{\{{\cal C}^*\}} 
 \rho^{\vert \CC^*\vert }\, { (n-1)(n-2) \cdots (n- \vert {\cal C}^*\vert )\over n^{\vert {\cal C}^*\vert }}
$$
and 
$$
\tilde  M_k^{(n,\rho)} \ge  \left( 1-{\rho\over n}\right)^k 
\left( 1- {4\rho\over n-\rho}\right)^k \tilde \CB_k(\rho)
= \tilde \CB_k(\rho)(1+O(\rho k/n))
$$
(cf. (5.13) and (5.15)). Then (3.5) follows.

%%%%%%%%%%%%%%%%%%%%%%%%%%%%%%%%%%%%%%%%%%%%%
\subsection{On the solution of  modified Lambert equation}

In this subsection we  repeat the reasoning  developed by N. G. de Bruijn \cite{DB} to study (2.12)
and prove the following statement.

\vskip 0.2cm
{\bf Lemma 5.4.} {\it Let $\tilde u(t)$ is the solution of Lambert-type equation
$\ti u(e^{\ti u} -1) = t$, $t>0$ (2.18). Then 
$$
\tilde u (t)=  \ln t - \ln \ln t + O \left( { \ln \ln t\over \ln t}\right), \quad t\to\infty.
\eqno (5.16)
$$

Proof.}
Omitting tildes,  we  rewrite equality 
$
u(e^u-1)=t
$
as 
$$
\ln \left( e^u - 1\right) = \ln t - \ln u.
\eqno (5.17)
$$
Assuming that $t>e^2$, we deduce from (5.17) that $u>1$. In the opposite case, $0<u\le 1$, we would get the upper bound
$\ln(e^u-1) \le \ln (e-1)< \ln 2$ that contradicts to (5.17).
Since $u>1$, then 
$
 \ln \left( e^u - 1\right) < \ln t 
$
and 
$$
0< \ln u< \ln\left( \ln t+1\right)
$$
and therefore
$$
\ln \left( e^u - 1\right)= \ln t + O\left( \ln \ln t\right), \quad t\to\infty.
$$
We denote $\ln t + O\left( \ln \ln t\right)= R$.
Then 
$$
u = \ln \left( e^R +1\right) = R + \ln \left( 1+ {1\over e^R}\right) = \ln t + O\left( \ln \ln t\right).
$$
Taking logarithms of  the both sides of this equality, we see that 
$$
\ln u = \ln \left( \ln t + O\left( \ln \ln t\right)\right)= \ln \ln t + O\left( { \ln \ln t/ \ln t}\right).
$$
Now we can deduce from (5.17) that 
$$
\ln \left( e ^u - 1\right) = \ln t - \ln \ln t + O \left( { \ln \ln t\over \ln t}\right), \quad t\to\infty
$$
and (5.16) follows. Lemma 5.4 is proved. 

%\subsection{
%Strongly restricted Bell numbers}
%%%%%%%%%%%%%%%%%%%%%%%%%%%%%%%%%%%%%%%%%%%%%%%%%%%%
%%%

%%%%%%%%%%%%%%%%%%%%%%

\end{document}